\documentclass[12pt]{amsart}
\usepackage{amssymb,amscd}
\usepackage[mathscr]{eucal}
\input epsf
\newtheorem{thm}{Theorem}[section]

\newtheorem{defin}[thm]{Definition}
\newtheorem{conjecture}[thm]{Conjecture}
\newtheorem{prop}[thm]{Proposition}

\newtheorem{cor}[thm]{Corollary}
\begin{document}

\newcommand{\C}{{\Bbb C}}
\newcommand{\HH}{{\Bbb H}}
\newcommand{\M}{{\mathcal M}}
\newcommand{\R}{{\Bbb R}}
\newcommand{\Z}{{\Bbb Z}}
\newcommand{\E}{{\mathcal E}}
\newcommand{\SU}{\o{SU}(2)}
\newcommand{\Hom}{\o{Hom}}
\newcommand{\Spin}{\o{Spin}(2)}
\newcommand{\Spinm}{\o{Spin}_-(2)}
\newcommand{\Pin}{\o{Pin}(2)}
\renewcommand{\O}{{\mathcal O}}
\renewcommand{\P}{{\Bbb P}}
\newcommand{\Q}{{\mathcal Q}}
\renewcommand{\o}{\operatorname}
\newcommand{\tr}{\o{tr}}
\title
{Topological Dynamics on Moduli Spaces, I}
\author{Joseph P. Previte \and Eugene Z. Xia}
\address{ School of Science,
Penn State Erie, The Behrend College,
Erie, PA 16563}
\address{Department of Mathematics,
University of Arizona, Tucson, AZ 85721}
\email{jpp@vortex.bd.psu.edu {\it (Previte)}, exia@math.arizona.edu {\it (Xia)}}
\date{\today} 
\maketitle
\begin{abstract}
Let $M$ be a one-holed torus with boundary  
$\partial M$ (a circle) and
$\Gamma$ the mapping class group of $M$ fixing $\partial M$. 
The group $\Gamma$ acts on ${\mathcal M}_{\mathcal C}(\SU)$ which
is the space of $\SU$-gauge equivalence classes of flat 
$\SU$-connections on $M$ with fixed holonomy on $\partial M$.
We study the topological dynamics of the $\Gamma$-action and 
give conditions for the individual $\Gamma$-orbits to be dense
in ${\mathcal M}_{\mathcal C}(\SU)$.  
\end{abstract}

\section{Introduction}
Let $M$ be a Riemann surface of genus $g$ with $m$
boundary components (circles).  Let 
$$
\{\gamma_1, \gamma_2,..., \gamma_m\} \subset \pi_1(M)
$$
be the elements in the fundamental group corresponding to
these $m$ boundary components.  
The space of $\SU$-gauge equivalence classes of 
$\SU$-connections $YM_2(\SU)$ is the well known Yang-Mills two space 
of quantum field theory.  Inside $YM_2(\SU)$ is the moduli space 
${\mathcal M}(\SU)$ of flat $\SU$-connections.

The moduli space ${\mathcal M}(\SU)$ has an interpretation that relates
to the representation space $\Hom(\pi_1(X),\SU)$ which is
an affine variety.  The group $\SU$ acts on $\Hom(\pi_1(M),\SU)$
by conjugation, and the resulting quotient space is precisely
$$
{\mathcal M}(\SU) = \Hom(\pi_1(M),\SU)/\SU.  
$$
Conceptually, the moduli space ${\mathcal M}(\SU)$ relates
to the semi-classical limit of $YM_2(\SU)$.

Assign each $\gamma_i$ a conjugacy class $C_i \subset \SU$ and let
$$
{\mathcal C} = \{C_1, C_2, ..., C_m \}.
$$
\begin{defin}~\label{def:1.1}
The relative character variety with respect to ${\mathcal C}$ is
$$
{\mathcal M}_{\mathcal C}(\SU) = \{[\rho] \in {\mathcal M} : \rho(\gamma_i) \in 
C_i, 1 \le i \le m \}.
$$
\end{defin}
The space ${\mathcal M}_{\mathcal C}(\SU)$ is compact, but possibly singular.
The set of smooth points of ${\mathcal M}_{\mathcal C}(\SU)$ possesses a natural
symplectic structure $\omega$ which gives rise 
to a finite measure $\mu$ on ${\mathcal M}_{\mathcal C}(\SU)$ (see \cite{Go1, Go2, Hu}).

Let  $\o{Diff}(M, \partial M)$ be the group of diffeomorphisms fixing
$\partial M$.  The mapping class  group $\Gamma$ is defined
to be $\pi_0(\o{Diff}(M, \partial M))$.  The group $\Gamma$
acts on $\pi_1(M)$ fixing the $\gamma_i$'s.
It is known that $\omega$ (hence $\mu$)
is invariant with respect to the $\Gamma$-action. In \cite{Go1} ,
 Goldman showed that with respect to the measure $\mu$,
\begin{thm} [Goldman] \label{thm:1}
The the mapping class group $\Gamma$ acts
ergodically on ${\mathcal M}_{\mathcal C}(\SU)$.
\end{thm}

Since ${\mathcal M}_{\mathcal C}(\SU)$ is a variety, one may also study the
topological dynamics of the mapping class group action.  
The topological-dynamical problem is considerably more delicate.  
To begin,
not all orbits are dense in ${\mathcal M}_{\mathcal C}(\SU)$.  If 
$\sigma \in \Hom(\pi_1(M),G)$
where $G$ is a proper closed subgroup of $\SU$ and $\tau \in \Gamma$, then
$\tau(\sigma) \in \Hom(\pi_1(M),G)$.   In other words, 
${\mathcal M}_{\mathcal C}(G)
\subset {\mathcal M}_{\mathcal C}(\SU)$ is invariant with respect to the 
$\Gamma$-action.
The main result of this paper is the following:
\begin{thm} \label{thm:2}
Suppose that $M$ is a torus with one boundary component and
$\sigma \in \Hom(\pi_1(M),\SU)$ such that $\sigma(\pi_1(M))$
is dense in $\SU$.
%is not contained in a group isomorphic to $C, D,$ or $\Pin$.
Then the $\Gamma$-orbit of the conjugacy class
$[\sigma] \in {\mathcal M}_{\mathcal C}(\SU)$
is dense in ${\mathcal M}_{\mathcal C}(\SU)$.
\end{thm}

The group $\SU$ is a double cover of $\o{SO}(3)$:
$$
p : \SU \longrightarrow \o{SO}(3).
$$
The group $\o{SO}(3)$ contains $\o{O}(2)$,
and the symmetry groups of the regular polyhedra: $T'$ (the tetrahedron),
$C'$ (the cube), and $D'$ (the dodecahedron).
Let $\Pin, T, C,$ and $D$ denote the groups 
$p^{-1}(\o{O}(2)),$ $ p^{-1}(T'),$ $p^{-1}(C'),$
and $p^{-1}(D')$, respectively.
The proper closed subgroups of $\SU$ consist of $T$, $C$, $D,$ and the 
closed subgroups of $\Pin$.  In particular, $T \subset C$
and $T \subset D$.  Suppose $\sigma \in \Hom(\pi_1(X),\SU)$.  Denote
$[\sigma]$ the corresponding $\SU$-conjugacy class in 
${\mathcal M}_{\mathcal C}(\SU)$.
Theorem~\ref{thm:2} implies that if $\sigma(\pi_1(M))$ is not contained in 
a group isomorphic to $C, D,$ or $\Pin$, then the $\Gamma$-orbit 
of the conjugacy class $[\sigma] \in {\mathcal M}_{\mathcal C}(\SU)$
is dense in ${\mathcal M}_{\mathcal C}(\SU)$.
%The main result of this paper is the following:
%\begin{thm} \label{thm:2}
%Suppose $M$ is a torus with one boundary component and 
%$\sigma \in \Hom(\pi_1(M),\SU)$ such that $\sigma(\pi_1(M))$
%is dense in $\SU$.
%%is not contained in a group isomorphic to $C, D,$ or $\Pin$.  
%Then the $\Gamma$-orbit of the conjugacy class 
%$[\sigma] \in {\mathcal M}_{\mathcal C}(\SU)$
%is dense in ${\mathcal M}_{\mathcal C}(\SU)$.
%\end{thm}

A conjugacy class in $\SU$ is determined by its trace.
For $M$ a torus with one boundary component, the moduli space 
${\mathcal M}(\SU)$
is a topological ball while ${\mathcal M}_{\mathcal C}(\SU)$ is
generically a smooth 2-sphere.  The mapping class group $\Gamma$
is generated by two Dehn twists $\tau_X, \tau_Y$.  With 
a proper change of coordinates, $\tau_X$ and $\tau_Y$ act on 
${\mathcal M}_{\mathcal C}(\SU)$ by rotations along two axes.  
The angle of rotation depends on the latitude of the circle 
along the respective axis.  The ergodicity theorem for 
${\mathcal M}_{\mathcal C}(\SU)$ follows 
since almost all the rotations are irrational multiples
of $2 \pi$.  However, in the context of topological dynamics, 
one must analyze the orbit of each class 
$[\rho] \in {\mathcal M}_{\mathcal C}(\SU)$
upon which one or both of the $\tau_X, \tau_Y$ actions are
rotations of rational multiples of $2 \pi$.

The proof of Theorem~\ref{thm:2} consists of two steps.  
The first is purely 
topological-dynamical in nature, concerning the case when the $\Gamma$-orbit
is infinite.  The second step deals with the cases where the $\Gamma$-orbits
are potentially finite and involves the theory of trigonometric 
Diophantine equations.
All in all, the proof is a delicate interplay of ideas
in geometric invariant theory \cite{Mu, Ne}, topological dynamics, and Diophantine
equations.  Incidentally, the proof also yields the well-known result
that the only proper closed subgroups $\SU$ are the closed subgroups
of $\Pin$ and the double covers of the automorphism groups of the
Platonic solids.

% Eug add in why one-holed is important
%
The following conjecture is the analogue of Theorem~\ref{thm:1}
in the category of topological dynamics.
Theorem~\ref{thm:2} is a major
stepping stone in the search of
a proof for this conjecture.
\begin{conjecture} \label{con:1.1}
Suppose that $M$ is a Riemann surface with boundary and
$\sigma \in \Hom(\pi_1(M),\SU)$ such that $\sigma(\pi_1(M))$
is dense in $\SU$.
%is not contained in a group isomorphic to $C, D,$ or $\Pin$.  
Then the 
$\Gamma$-orbit of the conjugacy class $[\sigma] \in {\mathcal M}_{\mathcal C}(\SU)$
is dense in ${\mathcal M}_{\mathcal C}(\SU)$.
\end{conjecture}

\centerline{\sc Acknowledgments}
 
We thank Professors Michael Brin, William Goldman, Larry Grove,
Kirti Joshi, David Levermore,
William McCallum, Michelle Previte, and Lawrence Washington
for insightful discussions during the course of this research.
Eugene Xia also thanks IH\'{E}S for providing
an excellent environment during the final phase of this 
research.
 
\section{Coordinates on the Moduli Space}
For the rest of this paper,
fix $M$ to be a torus with one boundary component.
We write $E$ for ${\mathcal M}(\SU)$
and $E_k$ for ${\mathcal M}_{\mathcal C}(\SU)$ such that $k = \tr(C)$,
where $C$ is the sole element in ${\mathcal C}$.
%and write ${\mathcal M}_k$ for ${\mathcal M}_{\mathcal C}(\SU)$ where $k = \tr K$
%and $K$ is the sole element in ${\mathcal C}$.
In this section, we briefly summarize some general properties
of $E$.  Consult \cite{Go1} for details.

The fundamental group $\pi_1(M)$ has a presentation
$$
\pi_1(M) = \langle X,Y,K | K = XYX^{-1}Y^{-1} \rangle
$$
where $K$ represents the element generated by the boundary component.
In particular, $\pi_1(M)$ is the free group generated by $X$ and $Y$.
Note
$$
E = \Hom(\pi_1(M), \SU)/\SU.
$$
The $\SU$-invariant polynomials \cite{Ne} on $\Hom(\pi_1(M), \SU)$ are
generated by the traces of the representations.  In particular,
a point $[\sigma] \in E$ is determined by
$$
x = \tr(\sigma(X)), y = \tr(\sigma(Y)), z = \tr(\sigma(X Y)).
$$
This provides a global coordinate chart:
$$
F : E \longmapsto \R^3
$$
$$
[\sigma] \stackrel{F}{\longmapsto} (\tr(\sigma(X)), 
\tr(\sigma(Y)), \tr(\sigma(X Y))).
$$
In addition, $k = \tr(\sigma(K))$ is given by the formula
\begin{equation} \label{eq0}
k = \tr(\sigma(K)) = x^2 + y^2 + z^2 - x y z - 2.
\end{equation}
The trace of every element in $\SU$ is in $[-2,2].$
In fact, one can show that
$$
E = \{(x,y,z) \in [-2,2]^3: -2 \le k \le 2\}.
$$
Let $f$ be the map
$$
f: E \longrightarrow [-2, 2]
$$
$$
f([\sigma]) = \tr(\sigma(K)).
$$
The fibre $f^{-1}(k)$ is precisely $E_k$ and is a smooth 2-sphere for 
each $-2 < k < 2$.
The fibre $f^{-1}(2)$ is a singular sphere while $f^{-1}(-2)$ consists of
one point. 

The mapping class group $\Gamma$ is generated by the maps 
$\tau_X$ and $\tau_Y$:
$$
\tau_X(X) = X \mbox{ and } \tau_X(Y) = YX
$$
$$
\tau_Y(X) = XY \mbox{ and } \tau_Y(Y) = Y.
$$
The induced action of $\Gamma$ on
$E$ preserves $E_k$.

With respect to the global coordinate, the actions of $\tau_X$ and
$\tau_Y$ can be described explicitly:
$$
\tau_X(x,y,z) = (x, z, xz - y)
$$
$$
\tau_Y(x,y,z) = (z, y, yz - x).
$$
The action of $\tau_X$ fixes $x$ and $k$, and preserves
the ellipse
$$
E_{x,k} = \{x\} \times \{ (y,z) : \frac{2+x}{4} (y+z)^2 + \frac{2-x}{4} (y-z)^2 = 2 + k 
- x^2 \}.
$$

The topological sphere $f^{-1}(1)$ is pictured below, decomposed into
ellipses.

\begin{figure}[h]
\centerline{\epsfysize=1.5in		\epsffile{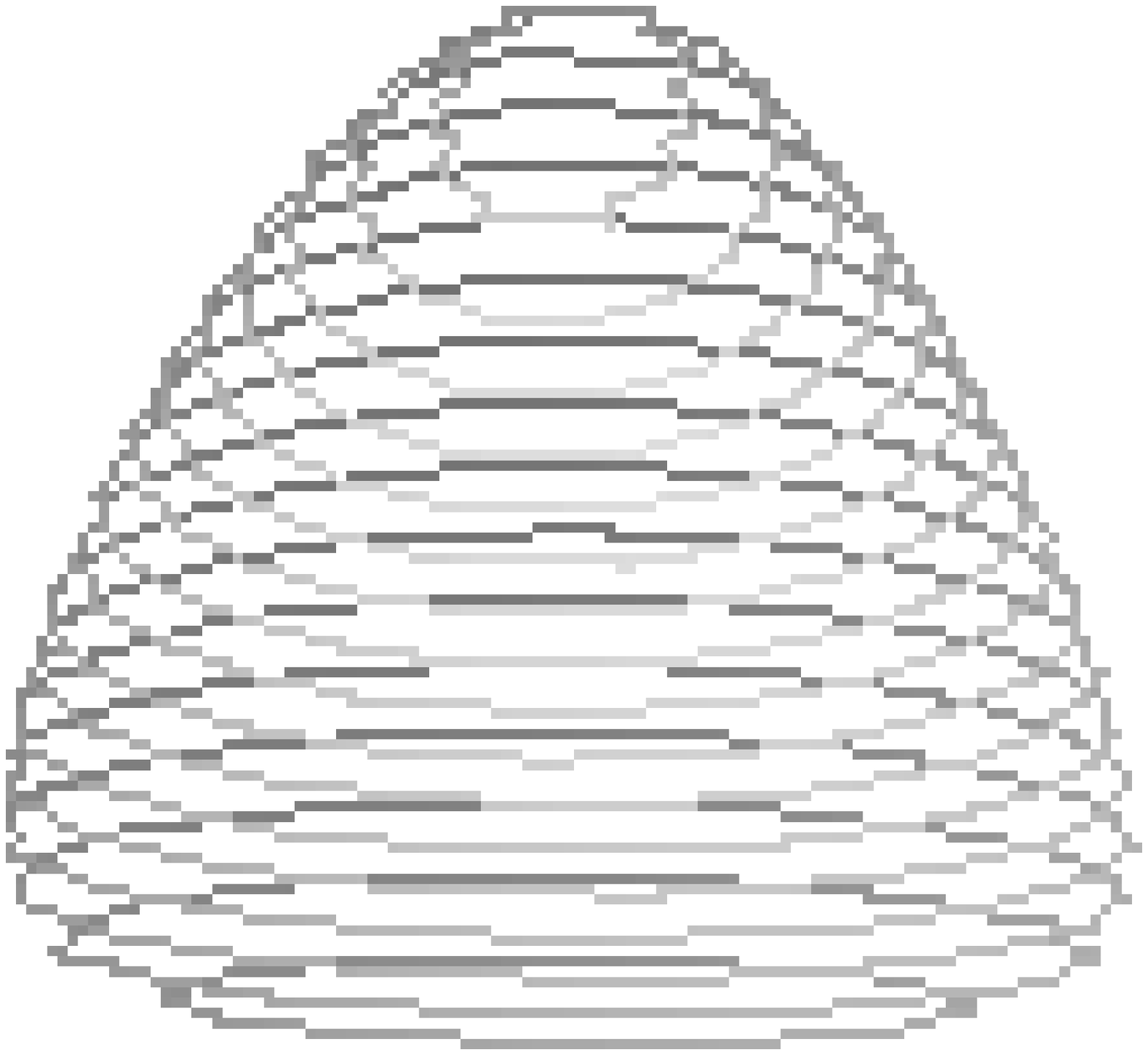}}
\caption{the topological sphere $E_1$}	
\end{figure}

The change of coordinates 
$$
\left\{
\begin{array}{lll}
\tilde{x} & = & x\\[2ex]
\tilde{y} & = & \frac{\sqrt{2-x} + \sqrt{2+x}}{2 \sqrt 2} y + 
\frac{\sqrt{2-x} - \sqrt{2+x}}{2 \sqrt 2}z \\[2ex]
\tilde{z} & = & \frac{\sqrt{2-x} - \sqrt{2+x}}{2 \sqrt 2} y + 
\frac{\sqrt{2-x} + \sqrt{2+x}}{2 \sqrt 2}z
\end{array}
\right.
$$
transforms $E_{x,k}$ into the circle
$$
E_{x,k} = \{\tilde{x}\} \times \{(\tilde{y},\tilde{z}) : 
\tilde{y}^2 + \tilde{z}^2 = 2 + k - \tilde{x}^2 \}.
$$
In this new coordinate system, $\tau_X$ acts as a rotation
by $\cos^{-1}(x/2)$.  In short, the sphere $E_k$ is
the union of circles
$$
E_k = \bigcup_{x} E_{x,k} ,
$$
and $\tau_X$ rotates (up to a coordinate transformation) 
each level set $E_{x,k}$ by an angle of 
$\cos^{-1}(x/2)$.

Under an analogous coordinate transformation, the action
$\tau_Y$ becomes a rotation of $E_{y,k}$ by an angle of
$\cos^{-1}(y/2).$

\section{The Closed Subgroup Representations}
All closed proper subgroups of $\SU$ are contained in $\Pin,$
$C,$ or $D$, where $C$ and $D$ are the double covers of the 
isometry groups of the cube and dodecahedron, respectively.
In this section, we classify the $\Pin$ representation classes
and produce a list of global coordinates for some
$C, D$ representation classes.  We shall prove later that
the list is complete up to some simplifying assumptions and
the variations allowed by
the following proposition:
\begin{prop} \label{prop:3.1}
Suppose $\sigma$ is a $G$ representation such that
$G \subset \SU$. If
$(x,y,z)$ are the global coordinates of $[\sigma]$,
then any permutation
of $(x,y,z)$ also corresponds to a $G$ representation class.
In addition, if $-I \in G$, then
the triples
$(-x,y,-z)$,
$(x,-y,-z)$,
$(-x,-y,z)$ also correspond to $G$ representation classes.
\end{prop}
\begin{proof}
Suppose $\sigma$ is a $G$ representation such that
$[\sigma]$ has global coordinates $(x,y,z)$.
Since $G$ is a group, the representation $\sigma'$, with
$$
\sigma'(X) = \sigma(XY), \ \ \ \sigma'(Y) = \sigma(X^{-1}),
$$
is also
a $G$ representation.  Moreover $[\sigma']$ has global coordinates 
$(z,x,y)$. The other permutations of coordinates are handled similarly.

Suppose $\sigma$ is a $G$ representation.  Since $-I \in G$, the
representation $\sigma'$, with
$$
\sigma'(X) = -\sigma(X), \ \ \ \sigma'(Y) = \sigma(Y),
$$
is also a $G$ representation and the global coordinates of $[\sigma']$ 
is $(-x, y, -z)$.  The other cases are handled similarly.
\end{proof}
 
If $-I \in G$, then
two classes
$$
[\sigma], [\sigma'] \in \Hom(\pi_1(M),G)/G
$$ 
are called
$S$-equivalent if their global coordinates differ from
one another as prescribed by Proposition~\ref{prop:3.1}.
If $-I$ is not in $G$, then two classes $[\sigma], [\sigma']$
are called $S$-equivalent if their global coordinates differ from
one another by a permutation.
 
\subsection{The $\Pin$ Representations}
The group $\Pin$ has two components, and
we write
$$
\Pin = \Spin \cup \Spinm,
$$
where $\Spin$ is the identity component of $\Pin$.
\begin{prop} \label{prop:3.2} A representation
$\sigma$ is a $\Spin$ representation if and only if
$\tr(\sigma(K)) = 2$.
\end{prop}
\begin{proof}
If $\sigma$ is a $\Spin$ representation, then
$$
\sigma(K) = \sigma(X Y X^{-1} Y^{-1}) = I,
$$
since $\Spin$ is abelian.  Hence,
$$
k = \tr(\sigma(K)) = \tr(I) = 2.
$$
If $k = \tr(\sigma(K)) = 2$, then
$$
I = \sigma(K) = \sigma(X) \sigma(Y) \sigma(X)^{-1} \sigma(Y)^{-1}.
$$
This implies that
the image of $\sigma$ is abelian, hence, is contained in $\Spin$.
\end{proof}
\begin{prop} \label{prop:3.3}
A representation $\sigma$ is a $\Pin$ representation and not
a $\Spin$ representation if and only if $k \neq 2$ and 
at least two of the three global
coordinates of $[\sigma]$ are zero.
\end{prop}
\begin{proof}
If $\sigma$ is a $\Pin$ representation and not
a $\Spin$ representation, then at least two of the following 
$$
\sigma(X), \sigma(Y), \sigma(XY)
$$
are in $\Spinm$.
Since $A \in \Spinm$ implies $\tr(A) = 0$, at least two of the
three global coordinates of $[\sigma]$ must be zero.
 
Suppose two of the three global
coordinates of $[\sigma]$ are zero, say $y, z = 0$. 
One easily finds a
$\Pin$ representation $\sigma'$ such that the global coordinates 
of $[\sigma']$ are $(x,0,0)$:
$$
\sigma'(X)=\begin{pmatrix}
\cos(\theta) & \sin(\theta) \\
-\sin(\theta) & \cos(\theta)
\end{pmatrix},
\sigma'(Y)=\begin{pmatrix}
0 & i \\
i & 0
\end{pmatrix},
$$
where $x = 2 \cos(\theta)$.  If $x \neq \pm 2$, then $k \neq 2$ and
$\sigma$ is not a $\Spin$ representation.
Since global coordinates are unique,
$$
[\sigma] =[\sigma'] \in \Hom(\pi_1(M),\SU)/\SU.
$$
Thus, $\sigma$ is a $\Pin$ representation
but not a $\Spin$ representation.
The proofs for the other cases are similar.
\end{proof}
 
This provides a complete characterization of the $\Pin$ 
representation classes.  
\begin{cor} \label{cor:3.2}
The space of $\Spin$ representation classes
consists precisely of $E_2$ (or ${\partial} E$). The
$\Pin$ representation classes consist of $E_2$ and the intersections of
the three axes with $E$.
For $-2 < k < 2$, there are exactly
six points corresponding to $\Pin$ representation classes
in $E_k$.
\end{cor}
 
\subsection{The $C$ and $D$ Representations}
Since $C$ (respectively, $D$) is finite, the $\Gamma$-orbit of a $C$ (respectively, $D$) representation class is finite.
One also notes that $-I$ is in $C$ and $D$.
 
We introduce the quaternionic model for $\SU$, namely
set $1,{\o i}, {\o j}, {\o k}$ as
$$
\left[
\begin{array}{cc}
1 & 0 \\
0 & 1 \\
\end{array}
\right], \ \
 \left[
\begin{array}{cc}
i & 0 \\
0 & -i\\
\end{array}
\right], \ \
\left[
\begin{array}{cc}
0 & 1\\
-1& 0 \\
\end{array}
\right], \ \
 \left[
\begin{array}{cc}
0 & i \\
i & 0\\
\end{array}
\right], \ \
$$
respectively. Then
$\SU = \{x+y{\o i}+z{\o j}+w{\o k} : x^2+y^2+z^2+w^2=1\},$
with the usual quaternionic multiplication.
 
For the rest of the paper, we fix the constants
$$
\left\{
\begin{array}{lll}
r & = & \frac{\sqrt{5} + 1}{4} \\[2ex]
s & = & \frac{\sqrt{5} - 1}{4}.
\end{array}
\right.
$$
 
Let
$$
\left\{
\begin{array}{lll}
T & = & \frac {\sqrt 2} 2 + \frac {\sqrt 2} 2{\o i} \\[2ex]
U & = & \frac {\sqrt 2} 2 + \frac {\sqrt 2} 2{\o j},
\end{array}
\right.
$$
$$
\left\{
\begin{array}{lll}
A & = & r+s {\o i}+\frac1 2 {\o k} \\[2ex]
B & = &-s +\frac 1 2  {\o i}-r{\o k}.
\end{array}
\right.
$$
%generates $D$.
Then
$$
\langle T, U \rangle \cong C \mbox{ and } \langle A, B \rangle \cong D.
$$
 
\begin{table}
\begin{center}
\begin{tabular}{|c|c|c|c|}
\hline
& & $S$-equivalence classes &\\
$X$ & $Y$ & $(x,y,z)$ & $k$\\
\hline
 
$T$ & $(TU)^{-1}$ & $(\sqrt 2,1,\sqrt 2)$ & $1$\\
 
$A^3B^8$ & $(AABA)^{-1}$ & $(1,1,1)$ & $0$\\
 
$B^{-1}$ & $AAB$ & $(-2s,2s,2s)$ &  $\frac{1 - \sqrt{5}}{2}$\\
 
$ABA^3B^2$ & $B$ & $(-2s,-2s,1)$ & $\frac{1 - \sqrt{5}}{2}$\\
 
$A^{-1}$ & $ABB$ & $(2r,-2r,-2r)$ & $\frac{1 + \sqrt{5}}{2}$\\
 
$ABAAB^6$ & $AABA$ & $(2r,1,2r)$ & $\frac{1 + \sqrt{5}}{2}$\\
 
$ABAA$ & $A$ & $(1,2r,1)$ & $1$\\
 
$AB^{-6}$ & $B^{-6}$ & $(1,2s,-1)$ & $1$\\
 
$ABAA$ & $A^{-1}$ & $(1,2r,2s)$ & $1$\\
\hline
 
\end{tabular}
\end{center}
\caption{Global coordinates for some $C$ and $D$ representation classes}
\end{table}

Table~1 is a list of
$S$-equivalence classes that come from
either $\langle T,U \rangle$
(cube) or  $\langle A,B \rangle$ (dodecahedron) representation classes.
We shall prove in Section~6 the following statement.
Suppose $(x,y,z)$ is a $C$ or $D$ representation class, but not a $\Pin$
representation class.  Then there exists $\gamma \in \Gamma$ such
that $\gamma(x,y,z)$ is in one of the $S$-equivalence classes in Table~1.
%$\tr(X)=x,$ $\tr(Y)=y$ and $\tr(XY)=z$.

\section{The Irrational Rotations and Infinite Orbits}

The Dehn twist $\tau_Y$ acts on the (transformed)
subsets $E_{y,k}$ via a rotation of angle  $\cos^{-1}(y/2)$.
The $y$-coordinates that yield finite orbits under $\tau_Y$ create a
filtration as follows:  
Let $Y_n \subset (-2, 2)$ such that
$y \in Y_n$ implies that if $(x,y,z)\in E$ are global coordinates
of a representation class, then the $\tau_Y$-orbit of $(x,y,z)$
is periodic with period greater than one but less than or equal 
to $n$.
This gives a filtration
$$
\emptyset = Y_2 \subset Y_3 \subset ... \subset Y_n \subset ...
$$
%and a grading:
%$$
%Gr(Y,n) = Y_n \setminus Y_{n-1}.
%$$
For example, $Y_2 = \emptyset$, $Y_3 = \{-1\}$, $Y_4 = \{-1, 0\}$,
etc. In particular, $Y_n$ is a finite set for every $n$.  
By symmetry, there exists a similar filtration $X_n,$ with $X_n=Y_n$ as sets.

Fix $-2 < k < 2$ and consider the two-dimensional sphere $E_k$.
The global coordinates provide an embedding of $E_k$ in $\R^3$ 
as a submanifold.  Hence $E_k$ inherits
a metric from the flat Riemannian metric on $\R^3$.  This provides
a distance function (metric) $d$ on $E_k$.  The metric $d$ generates
the usual topology on $E_k$.
Note that there are two points on 
$E_k$ that are fixed by $\tau_Y$.  These points correspond to
$\Pin$ representation classes.

\begin{defin}
For $\epsilon > 0$, a set $U$ is 
$\epsilon$-dense if for each $p \in U \subset E_k$, 
there exists a point $q \in U$ such that  
$
0 < d(p, q) < \epsilon.
$
\end{defin}
Let $\epsilon > 0$.
Since $E_k$ is compact,
there is an $M(\epsilon)$ such that
$n \ge M(\epsilon)$ implies that every $\tau_Y$-orbit $\O_y$ in $E_{y,k}$
is $\epsilon$-dense for any $y \not\in Y_{M(\epsilon)}$.
Let $N(\epsilon)$ be the cardinality of
$Y_{M(\epsilon)}$.
%Thus, there is a  finite number $N(\epsilon)$ of $y$-coordinates
%having periodic $\tau_Y$-orbits with the
%maximum distance (on $E_k$) separating
%points on a $\tau_Y$-orbit greater than 
%or equal to $\epsilon .$ 
%and a grading
%$Gr(X,n)=Gr(Y,n)$ with respect to $\tau_X$.

\begin{prop} \label{prop:4.1}
Let $(x_0,y_0,z_0) \in E_k.$ 
Suppose that one of $\cos^{-1}(x_0/2)$ or $\cos^{-1}(y_0/2)$ is 
an irrational multiple of $\pi$. Then the $\Gamma$-orbit
of $(x_0,y_0,z_0)$ is dense in $E_k$.
\end{prop}
\begin{proof}

Suppose that $\cos^{-1}(x_0/2)$ is an irrational multiple of $\pi$. 
Let $\epsilon > 0$
and  $(x_*,y_*, z_*) \in E_k$ which does not correspond to a
$\Pin$ representation class.
Since $\cos^{-1}(x_0/2)$ is an irrational multiple of $\pi$,
$\tau_X$ acts on the (transformed) subset of $E_{x_0,k}$,
by an irrational rotation.
By the compactness of $E_k$, there exists a 
$y$-value, $y_1\neq y_*$ and $\delta > 0$ such that 
$E_{y_*,k}$ is in the $\epsilon$-neighborhood of $E_{y_1,k}$ and
$
0 < \delta < d(E_{y_1,k}, E_{y_*,k}).
$
%$\epsilon$-neighborhood of $E_{y_1,k}$ covers $E_{y_*,k}.$
%In case $y_* \neq 0$, we also require that $|y_*| > |y_1|$.
We first consider the special case where
there exists an integer $J$ such that
the $y$-coordinate of $\tau_X^J(x_0,y_0,z_0)$ is strictly between 
$y_1$ and $y_*$.  Since $\cos^{-1}(x_0/2)$ is 
an irrational multiple of $\pi,$ there are infinitely many 
integers $J_i$ such that the $y$-coordinate of 
$\tau_X^{J_i}(x_0,y_0,z_0)$ is strictly between $y_1$ and $y_*$. 
Choose $J_i$ such that 
the $y$-coordinate of $(\tau_X)^{J_i}(x_0,y_0,z_0)$ 
is not in $Y_{M(\epsilon)}$.
%to correspond to a $y$-coordinate 
%$y_{M_i}$ so that the maximum distance separating
%points on a $\tau_Y$-orbit of $E_{y_{M_i},k}$ is less than $\epsilon.$
By the triangle inequality, there is some point on the 
$\tau_Y$-orbit of $(\tau_X)^{J_i}(x_0,y_0,z_0)$ that is
at most $2 \epsilon $ from $(x_*,y_*,z_*).$

We now prove the proposition in general.
%The following argument covers the preceding case as well as the case
%where there does not exist an integer $J$ such that
%the $y$-coordinate of $\tau_X^J(x_0,y_0,z_0)$ is strictly between 
%$y_1$ and $y_*.$
%Essentially, this means that
%there are no solutions to the inequality 
%$$
%|(x_0)^2 + (y_*)^2 + z^2-x_0y_*z-2-k|< \epsilon 
%$$ 
%for any 
%$(x_0,y_*,z) \in E_{x_0,k}.$  
Since $(x_*,y_*,z_*)$
satisfies equation (\ref{eq0}) and 
is not a $\Pin$ representation class, $E_{y_*,k}$ is a circle.
Hence, we must have that $(0,y_*,z') \in E_{y_*,k}$
for some $z'\neq 0$.
Therefore, there exists $x_2\neq 0$ such that 
$E_{x_2,k}$ intersects $E_{y_*,k}.$

Choose $\epsilon' = \frac{\delta}{N(\epsilon)+2}$.
By the filtration $X_n,$ 
the set $X_{M(\epsilon')}$ contains
all $x$-values that have the following properties:
\begin{enumerate}
\item $E_{x,k}$ intersects $E_{y_*,k}$
\item There is a point in $E_{x,k}$
whose  $\tau_X$-orbit has at most $N(\epsilon)$ points 
with distinct $y$-coordinates between $y_1$ and $y_*$.
\end{enumerate}

Note that the $x$-coordinate
of $\tau_Y(\tau_X)^J(x_0,y_0,z_0)$ is the
$z$-coordinate of $(\tau_X)^J(x_0,y_0,z_0).$
Since $\cos^{-1}(x_0/2)$ is an irrational multiple of $\pi,$
there is an infinite sequence of numbers $J_i$ such that
$|x_{J_i}| < |x_2|$, where $x_{J_i}$ is the the $x$-coordinate
of $\tau_Y (\tau_X)^{J_i}(x_0,y_0,z_0)$.
This forces
$E_{x_{J_i},k}$ to intersect $E_{y_*,k}$.
Of these, choose $J$ such that $x_{J}$
is not in $X_{M(\epsilon')}.$ Thus, the 
 $\tau_X$-orbit 
of $\tau_Y(\tau_X)^J(x_0,y_0,z_0)$ 
has at least $N(\epsilon)+1$ points 
with distinct $y$-coordinates between $y_1$ and $y_*$.

Now at most $N(\epsilon)$ values of $y$
yield $\tau_Y$-orbits that are not
$\epsilon$-dense. Thus, there exists a 
point $(\hat x, \hat y, \hat z)$ on the $\tau_X$-orbit
of $\tau_Y(\tau_X)^J(x_0,y_0,z_0)$ 
such that $\hat{y}$ is between $y_1$ and $y_*$, moreover,
the $\tau_Y$-orbit of $(\hat x, \hat y, \hat z)$ is
$\epsilon$-dense.
%$$
%d(E_{\hat y,k}, E_{y_*,k}) < \epsilon.
%$$
Since the $\epsilon$-neighborhood of $E_{\hat y,k}$ covers
$E_{y_*,k}$, 
 some point in the $\tau_Y$-orbit of $(\hat x, \hat y, \hat z)$
comes within $2 \epsilon$ of $(x_*,y_*,z_*).$ 
Finally, by Corollary~\ref{cor:3.2}, the set of $\Pin$ representation classes
in $E_k$ consists of six discrete points.  This implies there is no
loss of generality in assuming that $(x_*,y_*,z_*)$ 
does not correspond to a $\Pin$ 
representation class.
A symmetric argument holds if $\cos^{-1}(y_0/2)$ is an 
irrational multiple of $\pi.$ 
\end{proof}

\begin{prop} \label{prop:4.2}
Suppose 
the $\Gamma$-orbit of $(x_0,y_0,z_0) \in E_k$
is infinite.
Then the $\Gamma$-orbit
of $(x_0,y_0,z_0)$ is dense in $E_k$.
\end{prop}
\begin{proof}
Let $(x_0,y_0,z_0) \in E_k$ have infinite $\Gamma$-orbit and 
 $(x_*,y_*, z_*) \in E_k$ which does not correspond to a
$\Pin$ representation class. 

There are two cases.
One possibility is that the $\Gamma$-orbit $\mathcal O$ 
has an infinite number of points on some
circle  $E_{y,k}$ (respectively, $E_{x,k}$).
Hence,  there is an infinite number of points on 
$\mathcal O \cap E_{y,k}$ that have distinct
$x$-coordinates. However, a priori, these points may not be dense in
 $E_{y,k}.$ One uses the infinite number of points on 
$\mathcal O \cap E_{y,k}$ with distinct
$x$-coordinates as in the proof of Proposition~\ref{prop:4.1}.

%The proposition then follows by an argument similar to that of
%Proposition~\ref{prop:4.1}.
%of Proposition~\ref{prop:4.1} now gives the desired result.
%One now uses the
%filtration $X_n$ as in the
%previous proposition applied to 
%one of these points where $\epsilon$ is 
%suitably small.

The other possible case is that no circle $E_{x,k}$ (or $E_{y,k}$)
has an infinite number of points on $\mathcal O.$ 
As in the previous case, there are an
infinite number of points on the
 $\Gamma$-orbit of $(x_0,y_0,z_0)$ 
having distinct $x$-coordinates. The proof again follows
similarly to the proof of Proposition~\ref{prop:4.1}.

\end{proof}

\section{Trigonometric Diophantine Equations}
It remains for us to classify the finite $\Gamma$-orbits.  These orbits
exist and can be constructed by taking $G$ representation classes 
where $G$ is a closed proper subgroup of $\SU$ as described in Section~3. 

The problem amounts to considering cases where
the rotations generated by $\tau_X$ and $\tau_Y$ are both rational.  
In such cases, an additional iteration is made as follows.  
By assumption, both
$\cos^{-1}(x/2)$ and $\cos^{-1}(y/2)$ are rational multiples
of $\pi$. Also, $ \cos^{-1}(z/2)$ is a rational multiple
of $\pi$ since the $x$-coordinate of $\tau_Y(x,y,z)$ is $z$. 
Since $\tau_X(x,y,z)=(x,z,xz-y)$,  
in order for the orbit to be finite,
$\cos^{-1}(\frac{xz-y}2)$ must be a rational multiple of $\pi.$
In particular, $x = 2 \cos (\theta_x),$
$y = 2 \cos (\theta_y),$ $z = 2 \cos (\theta_z)$
and $xz-y = 2 \cos (\theta_{xz-y}),$
where all angles are rational multiples of $\pi$ in $[0, \pi].$ 
Hence,
$$
2 \cos (\theta_x) \cos (\theta_z)-  \cos (\theta_y) = \cos(\theta_{xz-y}),
$$
or
\begin{eqnarray} \label{eq1}
\cos (\theta_x+\theta_z) +\cos (\theta_x-\theta_z)-  
\cos (\theta_y) = \cos(\theta_{xz-y}),
\end{eqnarray}
where all angles are rational multiples of $\pi$,
$0\le \theta_x+\theta_z \le 2\pi$,
and $-\pi \le \theta_x-\theta_z \le \pi.$
Similarly, the action of $\tau_Y$ gives
\begin{eqnarray} \label{eq3}
\cos (\theta_y+\theta_z) +\cos (\theta_y-\theta_z)-  
\cos (\theta_x) = \cos(\theta_{yz-x}),
\end{eqnarray}
where all angles are rational multiples of $\pi$,
$0\le \theta_y+\theta_z \le 2\pi$,
and $-\pi \le \theta_y-\theta_z \le \pi.$
Equations (\ref{eq1}) and (\ref{eq3}) are referred to as the $\tau_X$-equation
and $\tau_Y$-equation at $(x,y,z)$ (or at $(\theta_x, \theta_y, \theta_z)$),
respectively.

\begin{prop} \label{prop:trivtrig}
Let $(x,y,z) \in E_k$ with $x,y,z$ all non-zero.
Suppose that two terms appearing in
equation (\ref{eq1}) cancel one another.
Then $k=2.$
\end{prop}
\begin{proof}
%We must have that  $\cos(\theta_y) = - \cos(\theta_{xz-y})$ or
%$\cos(\theta_y) = \cos(\theta_x \pm \theta_z).$
Suppose $\cos(\theta_y) = - \cos(\theta_{xz-y}).$
Then $xz-y = -y$ which
implies that $x=0$ or $z=0,$ a contradiction to the
assumption that $x,y,z$ are all non-zero.

Suppose $\cos(\theta_y) = \cos(\theta_x + \theta_z)$.
Recall that $(x,y,z)$ satisfies equation (\ref{eq0}).
Thus,
\begin{eqnarray*}
k & = & 4\cos(\theta_x)^2 + 4\cos(\theta_x+\theta_z)^2+4\cos(\theta_z)^2\\
& & -8\cos(\theta_x)\cos(\theta_x+\theta_z)\cos(\theta_z)-2\\
%& = & 4\cos(\theta_x)^2 + 4\cos(\theta_x+\theta_z)^2+4
%\cos(\theta_z)^2\\
%& & -4[\cos(\theta_x+\theta_z)+\cos(\theta_x-\theta_z)]\cos(\theta_x+\theta_z)
%-2\\
& = & 4\cos(\theta_x)^2+4\cos(\theta_z)^2-4\cos(\theta_x-\theta_z)
\cos(\theta_x+\theta_z)\\
& = & 2\cos(2\theta_x)+2+2\cos(2\theta_z)-2\cos(2\theta_x)-2\cos(2\theta_z)=2.
\end{eqnarray*}

A similar argument applies in the case
$\cos(\theta_y)= \cos(\theta_x - \theta_z).$
\end{proof}

A symmetric argument shows that if two terms appearing in
equation (\ref{eq3}) cancel one another, then $k=2.$
For the remainder of this paper, we assume that $(x,y,z) \in E$ does not
correspond to a $\Pin$ representation class.
This implies that $k\neq 2$.  In addition,
we may assume that
all coordinates of $(x,y,z)\in E_k$ are non-zero: For
if  $x=0,$ then $y,z$ must both be non-zero by
Proposition~\ref{prop:3.3}. The point
$\tau_Y(0,y,z) =(z,y,yz)$ has all non-zero entries.

Equation (\ref{eq1})
is an at most four-term Diophantine equation, the solutions to which
are few as shown by Conway and Jones.

\begin{thm} [Conway, Jones] \cite{Co} \label{thm:3}
Suppose that we have at most four distinct rational multiples of $\pi$
lying strictly between $0$ and $\pi/2$ for which some linear
combination of their cosines is rational, but no proper subset
has this property. That is,
$$
A \cos( a) + B \cos(b) + C \cos( c) + D \cos(d) = E,
$$
for $A, B, C, D, E$ rational and $a,b,c,d \in (0, \pi/2)$
rational multiples of $\pi$.
Then the appropriate linear combination is proportional
to one from the following list:
\begin{eqnarray*}
\cos(\pi/3) & = & 1/2\\
\cos(t+\pi/3)+\cos(\pi/3-t)-\cos(t) & = & 0  \ \ (0< t <\pi/6)\\
\cos(\pi/5) - \cos(2\pi/5) & = & 1/2\\
\cos(\pi/7) - \cos(2\pi/7) + \cos(3\pi/7) & = & 1/2\\
\cos(\pi/5) - \cos(\pi/15) + \cos(4\pi/15) & = & 1/2\\
-\cos(2\pi/5)+\cos(2\pi/15)-\cos(7\pi/15) & = &  1/2\\
\cos(\pi/7) + \cos(3\pi/7) - \cos(\pi/21) +\cos(8\pi/21) & = & 1/2\\
\cos(\pi/7) - \cos(2\pi/7) + \cos(2\pi/21) -\cos(5\pi/21) & = & 1/2\\
-\cos(2\pi/7) + \cos(3\pi/7) + \cos(4\pi/21) +\cos(10\pi/21) & = & 1/2\\
-\cos(\pi/15) + \cos(2\pi/15) + \cos(4\pi/15) -\cos(7\pi/15) & = & 1/2.\\
\end{eqnarray*}
\end{thm}
 
The non-zero cosine terms in equation (\ref{eq1}) are not necessarily 
in $(0,\pi/2)$.
By applying the identities $\cos(\pi/2-t)= - \cos(\pi/2+t)$
and $\cos(\pi-t)=\cos(\pi+t)$, we derive from 
equation (\ref{eq1}) a new four-term cosine equation whose
arguments are in $[0,\pi/2]$.
That is, by a possible change of sign,
each term in equation (\ref{eq1}) may be rewritten with angle
in $[0,\pi/2].$   

Equation (\ref{eq1})
cannot correspond to the last four equations appearing in Theorem~\ref{thm:3}
since these equations have five non-zero terms.
By Proposition \ref{prop:trivtrig}, we may assume that
the angles appearing in equation (\ref{eq1})
are all distinct.
For if two or more angles are the same, then after combining
terms, the resulting equation must be proportional 
to the first equation in Theorem \ref{thm:3}.
This leads to $k=2$ by Proposition \ref{prop:trivtrig}.

\begin{prop} \label{prop:5.2}
Suppose $(x,y,z) \in E_k$ are not the global coordinates of a $\Pin$
representation class and 
with $x,y,z$ all non-zero.
 Suppose that
the $\Gamma$-orbit of $(x,y,z)$ is finite and
that
some angle in equation (\ref{eq1}) or equation (\ref{eq3})
is an integer multiple of $\pi$. Then
$(x,y,z)$ is S-equivalent to a triple appearing in 
Table ~1.
\end{prop}

\begin{proof}
%We now analyze the cases in which 
%some angle in equation (\ref{eq1})
%is an integer multiple of $\pi$.
%cannot be rewritten as an equation
%which satisfies the hypotheses of Theorem \ref{thm:3}.
By the assumption $k\neq 2$,
the only way that equation (\ref{eq1}) can have angles
equal to an integer multiple of $\pi$ is if $\theta_x-\theta_z =0$
or $\theta_x+\theta_z=\pi,$
i.e. $x=z$ or $x=-z$. Note that both cannot happen simultaneously.

Suppose $\theta_x+\theta_z=\pi$.  Then 
equation (\ref{eq1}) becomes
\begin{equation} \label{eq5}
\cos (\theta_x-\theta_z) -  \cos (\theta_y) - \cos(\theta_{xz-y})=1.
\end{equation}
Theorem~\ref{thm:3}  and the
assumption $y\neq 0$ lead to the following cases:
 
(A) Two terms in equation (\ref{eq5}) correspond to the first
equation of Theorem \ref{thm:3}, with one term
equal to zero (angle $\pi/2$).

Since $y \neq 0,  $ we must have that $\cos (\theta_y)=-\frac 1 2,$
so  $y = -1.$ 
Now either $\cos(\theta_{xz-y}) = 0$ or $-\frac 1 2$.
If $\cos(\theta_{xz-y}) = 0,$ then 
$xz = -1,$ hence $x^2 = 1$ which yields the triples
$(1,-1,-1)$ and $(-1,-1,1).$
If $\cos(\theta_{xz-y}) = \frac 1 2,$
then $xz = -2,$ or $x = \pm \sqrt 2.$
The resulting  triples $(x,y,z)$ are
$(\sqrt 2,-1,-\sqrt 2)$ 
and $(-\sqrt 2,-1,\sqrt 2).$
Note that all of the above triples belong to  an $S$-equivalence class
appearing in Table~1.

(B) Two terms in equation (\ref{eq5}) correspond to the third equation
 of Theorem
\ref{thm:3} while the remaining term corresponds to the first equation in
 Theorem \ref{thm:3}   .
The resulting triples are:
$$
(-2s, 2s, 2s), 
(2s, 2s, -2s),
(-2s,-1,2s),
(2s,-1,-2s),
$$
$$
(-2r,-2r,2r),
(2r, -2r, -2r),
(-2r,-1,2r),
(2r, -1,-2r),
$$
$$
(-1, -2r,1),
(1,-2r,-1),
(-1,2s, 1),
(1, 2s, -1).
$$
Note that all above triples belong to  an $S$-equivalence class
appearing in Table~1. 
% Later we will show that each of these triples correspond to $C$ or $D$ subrepresentations

Suppose $\theta_x-\theta_z =0$.  Then 
equation (\ref{eq1}) becomes
\begin{equation} \label{eq6}
\cos (\theta_x+\theta_z) -  \cos (\theta_y) - \cos(\theta_{xz-y})=-1.
\end{equation}
Theorem~\ref{thm:3} and the
assumption $y\neq 0$ lead to the following cases.

(A) Two terms in equation (\ref{eq6}) correspond to the first
equation of Theorem~\ref{thm:3}, with one term
equal to zero (angle $\pi/2$).
The various possibilities lead to the triples:
$$
(1,1,1), (-1,1,-1),(\sqrt 2,1,\sqrt 2),
(-\sqrt 2,1,-\sqrt 2).
$$
 
(B)  Two terms in equation (\ref{eq6}) correspond to 
the third equation of Theorem
\ref{thm:3} while the other corresponds to the first equation.
 The various possibilities lead to the triples:
$$
(2r, 2r, 2r),
(-2r, 2r, -2r),
(2r,1,2r),
(-2r, 1, -2r),
$$
$$
(1,2r,1),
(-1,2r,-1),
(1,-2s,1),
(-1,-2s,-1),
$$
$$
(2s,-2s,2s),
(-2s,-2s,-2s),
(2s,1,2s),
(-2s,1,-2s).
$$
Again, the $S$-equivalence classes of these triples appear in Table~1.
% Again, we will show that each of these triples correspond to $C$ or $D$ subrepresentations  
A similar argument holds if 
some angle in  equation (\ref{eq3})
is an integer multiple of $\pi$, i.e. $y = \pm z.$
\end{proof}

Henceforth, we assume that all angles in non-zero cosine terms
appearing in equation (\ref{eq1}) are distinct and
not integer multiples of $\pi$. Under these assumptions,
equation (\ref{eq1}) can be rewritten as an equation that satisfies
the hypotheses of the following which is a special case of Theorem~\ref{thm:3}:

\begin{thm} [Conway, Jones] \cite{Co} \label{thm:4}
Suppose that we have at most four distinct rational multiples of $\pi$
lying strictly between $0$ and $\pi/2$ for which some linear
combination of their cosines is zero, but no proper subset
has this property. That is,
$$
A \cos(a) + B \cos(b) + C \cos( c) + D \cos(d) = 0,
$$
for $A, B, C, D$ rational and $a,b,c,d \in (0, \pi/2)$ 
rational multiples of $\pi$.
Then the linear combination is
proportional to one from the following list:
 
\begin{eqnarray*}
\cos(t+\pi/3)+\cos(\pi/3-t)-\cos(t) & = & 0  \ \ (0< t <\pi/6)\\
\cos(\pi/5) - \cos(2\pi/5) - \cos(\pi/3) & = & 0\\
\cos(\pi/7) - \cos(2\pi/7) + \cos(3\pi/7) - \cos(\pi/3) & = & 0\\
\cos(\pi/5) - \cos(\pi/15) + \cos(4\pi/15) -\cos(\pi/3) & = & 0\\
-\cos(2\pi/5)+\cos(2\pi/15)-\cos(7\pi/15)-\cos(\pi/3) & = & 0\\
\end{eqnarray*}
\end{thm}

\section{Proof of Theorem~\ref{thm:2}}
In this section we prove the following proposition which
in turn proves Theorem~\ref{thm:2}:
\begin{prop} \label{prop:e3}
Let $(x,y,z)\in E_k$ with $x,y,z$ non-zero.
%Suppose that  equation (\ref{eq1}) corresponds to one of 
%the equations
%of Theorem \ref{thm:4}.
Then the $\Gamma$-orbit of $(x,y,z)$ is infinite
or there is $\gamma \in \Gamma$ 
such that  $\gamma(x,y,z)$ is $S$-equivalent to a triple in Table ~1.
\end{prop}
An immediate consequence of Proposition~\ref{prop:e3} is,
\begin{cor}
Suppose $(x,y,z)$ is a $C$ or $D$ representation class, but not a $\Pin$
representation class.  Then there exists $\gamma \in \Gamma$ such
that $\gamma(x,y,z)$ is in one of the $S$-equivalence classes in Table~1.
\end{cor}

The proof of Proposition~\ref{prop:e3} presented here is lengthy and
highly computational.  We begin by outlining
the overall strategy.  Consider all triples $(x,y,z)$ 
that arise from the solutions
of equation (\ref{eq1}) provided by Theorem~\ref{thm:4}.  For a
triple $(x,y,z)$ to have a finite $\Gamma$-orbit,
the four-term trigonometric equations that arise from repeated
applications of $\tau_X$ or $\tau_Y$ must have
solutions provided by Theorem~\ref{thm:4} or violate the hypotheses
of Theorem ~\ref{thm:4}.
We prove Proposition~\ref{prop:e3}
by showing that all triples $(x,y,z)$ that arise from the 
 solutions
of equation (\ref{eq1}) provided by Theorem~\ref{thm:4}
fall into one of the
following three categories:
\begin{enumerate}
\item $(x,y,z)$ has $k \ge 2$ or one of the global coordinates
$(x,y,z)$ is zero.
\item $(x,y,z)$ belongs to one of the $S$-equivalence classes
in Table~1. Hence, corresponds to a $C$ or a $D$ representation class.
\item $(x,y,z)$ has infinite $\Gamma$-orbit with $-2 < k < 2$. Hence, 
has dense $\Gamma$-orbit by Proposition~\ref{prop:4.2}.
\end{enumerate}

%We now examine all
%representation classes that arise from the solutions
%to equation (\ref{eq1}) as
%provided by Theorem~\ref{thm:4}. 

\begin{defin}
The pairs of angles $(\theta_z, \theta_x)$,
$(\pi - \theta_x, \pi - \theta_z)$, and
$(\pi - \theta_z, \pi - \theta_x)$ are called the symmetric,
dual, and dual-symmetric pairs of $(\theta_x, \theta_z)$, respectively.
\end{defin}

Recall that by a possible change of sign,
each term in equation (\ref{eq1}) may be rewritten with angle
in $[0,\pi/2].$   Therefore, for fixed
$a,b \in [0, \pi/2]$, we obtain the following eight systems of equations
$$
\left\{
\begin{array}{lll}
\cos(a) & = & \pm \cos(\theta_x+\theta_z)\\[2ex]
\cos(b) & = & \pm \cos(\theta_x-\theta_z)
\end{array}
\right.
\ \ \ \
\left\{
\begin{array}{lll}
\cos(b) & = & \pm \cos(\theta_x+\theta_z)\\[2ex]
\cos(a) & = & \pm \cos(\theta_x-\theta_z),
\end{array}
\right.
$$
for $\theta_x,\theta_z \in (0,\pi).$
The following, together with their dual, symmetric, and dual-symmetric
pairs, are all possible pairs $(\theta_x,\theta_z)$ 
that satisfy
one of the above eight systems of equations:
\begin{center}
$
(\frac{a+b}2,|\frac{a-b}2|),
(\pi-\frac{a+b}2,|\frac{a-b}2|),
$

$
(\pi/2-\frac{a+b}2, \pi/2+\frac{a-b}2),
(\pi/2-\frac{a+b}2, \pi/2+\frac{b-a}2).
$

\end{center}

% For the rest of the section, all equation numbers refer
% to the equations in Theorem \ref{thm:4}.  
We prove in detail the cases in which equation (\ref{eq1}) corresponds to
equation 2 or 3 in Theorem \ref{thm:4}.  The argument for equation 3 
is the simplest and exemplifies the primary techniques used in the other cases.
The case of equation 2 involves
a free parameter $t$, hence, is somewhat more involved than the 
others.  For the other cases, we simply enumerate all the possible
solutions and categorize them according to the above mentioned categories.

We first consider the case of 
$$
\cos(\pi/7) - \cos(2\pi/7) + \cos(3\pi/7) - \cos(\pi/3)  =  0
$$
in detail.
Consider $a =\pi/7$ and
$b =2\pi/7.$ As given above,  
the possibilities for $(\theta_x,\theta_z)$ are
$$
(3\pi/14, \pi/14), (11\pi/14, \pi/14), (2\pi/7, 3\pi/7),(2\pi/7, 4\pi/7),
$$
along with their  symmetric, dual, and dual-symmetric pairs. 
%The equation obtained by applying  $\tau_Y$ to $(x,y,z)$ is:
%$$
%\cos (\theta_y+\theta_z) +\cos (\theta_y-\theta_z)-  \cos (\theta_x) =
%\cos(\theta_{yz-x}).
%$$
Suppose first that $\theta_y = \pi/3$ (respectively,  $2\pi/3$).
Then the angles $\theta_y\pm \theta_z$  are rationally related to
$\pi$ by reduced rationals with denominators
$21$ or $42$. Thus, for $(x,y,z)$ to have a finite
$\Gamma$-orbit,
the $\tau_Y$-equation 
at $(\theta_x, \theta_y, \theta_z)$ must correspond to 
equation 1.
Since $x \neq 0$
this equation can be rewritten as equation 1 
only if $yz-x=0,$ or $y = \frac x z.$
Then
$y \neq \frac x z$ for  $\theta_y= \pi/3$ (or $2\pi/3),$ 
and
$x,z$
as given by the
four pairs listed above 
(as well as their dual, symmetric, and dual-symmetric pairs).
Thus the $\tau_Y$-equation 
at $(\theta_x, \theta_y, \theta_z)$ must violate the
hypotheses of Theorem~\ref{thm:4}.
 
If  $\theta_{xz-y}=\pi/3$ or $2\pi/3$
then, consider the point $\tau_X(\theta_x,\theta_y,\theta_z)
=(\theta_x,\theta_z,\theta_{xz-y})$ (in angle notation).
The $\tau_Y$-equation at $(\theta_x,\theta_z,\theta_{xz-y})$ is
$$
\cos (\theta_z+\theta_{xz-y}) +\cos (\theta_z-\theta_{xz-y})-  \cos (\theta_x) =
\cos(\theta_{z(xz-y)-x}),
$$
but here $xz-y$ is playing the role of $y$ 
in the argument given above. 

This type of argument works for 
the other five cases (e.g. $a=\pi/7,\ \ b=3\pi/7$, etc.). 

Now we cover the case
$$
\cos(t+\pi/3)+\cos(\pi/3-t)-\cos(t) = 0 \ \ \ (0 < t < \pi/6).
$$

Recall that we are under the standing assumption $x,y,z \neq 0.$
There are six cases.
 
Case 1: For $a=t+\pi/3,$ and  $b=\pi/3-t,$
the possibilities are: 
$$
(\pi/3,t), (2\pi/3, t),
(\pi/6,\pi/2+t),(\pi/6,\pi/2-t),
$$ 
and their dual, symmetric, and
dual-symmetric pairs. Since $y\neq 0$,
$\theta_y$ corresponds to $t$ and $xz-y=0.$ 

The first two pairs, together with their duals
$(\pi/3,\pi-t), (2\pi/3,\pi- t),$
yield $y=\pm z$.  By Proposition \ref{prop:5.2},  these triples
are either $S$-equivalent to those appearing in
Table 1 or have an  infinite $\Gamma$-orbits.

In the above symmetric (dual-symmetric) 
pairs, we have $y=\pm x$ and $z=1.$ 
The $\tau_X$ preimage
of $(x,\pm x, 1)$ is $(x, \pm x^2 - 1, \pm x)$ 
(signs taken together). Since
$ \pm x^2 - 1 \neq 0,$ 
we may apply Proposition \ref{prop:5.2}
to this triple as above. The same holds
for $z=-1.$

Now consider the last two pairs. 
Note that the angle $\theta_x=\pi/6$ does not appear in equations 1-5.
Therefore, any triple associated with these pairs will have an
infinite $\Gamma-$orbit. 
For the symmetric pairs, i.e. $\theta_z=\pi/6$ (respectively, $5\pi/6$),  
note that
$\tau_X(\theta_x,\theta_y,\theta_z)= (\theta_x,\theta_z,
\theta_{xz-y}).$ The $\tau_Y$-equation at $(\theta_x,\theta_z,\theta_{xz-y})$
cannot be put into the form of  equations
1-5, as above. This argument also applies to the duals of all such
pairs. 

%add in that cant be eq2

Case 2:  For $a=t+\pi/3,$ and  $b=t,$
the possibilities are: 
$$
(t+\pi/6,\pi/6), (5\pi/6-t, \pi/6),
(\pi/3-t,2\pi/3),(\pi/3-t,\pi/3),
$$ 
and their dual, symmetric, and
dual-symmetric pairs. 
Since $y\neq 0$,
$\theta_y$ corresponds to $\pi/3-t$ and $xz-y=0.$

As in the previous case, 
the last two pairs together with their dual, symmetric, and 
dual-symmetric pairs
yield triples (or triples whose $\tau_X$ preimage) are
either $S$-equivalent to those appearing in
Table 1 or have infinite $\Gamma$-orbits.
The argument for the first two pairs is also
similar to that in the previous case.

Case 3:  For $a=\pi/3-t,$ and  $b=t,$
the possibilities are: 
$$
(\pi/6, \pi/6-t), (5\pi/6, \pi/6-t),
(\pi/3, 2\pi/3-t),(\pi/3, t+\pi/3),
$$ 
together with their dual, symmetric, and
dual-symmetric pairs. 
An argument similar to the one given in case 1 holds.

Case 4: For $a=\pi/3+t,$ and $b= \pi/2,$
the possibilities are: 
\begin{center}
$
(5\pi/12+t/2, \pi/12-t/2), (7\pi/12-t/2, \pi/12-t/2),
$
$
(\pi/12-t/2,5\pi/12+t/2),(\pi/12-t/2,7\pi/12-t/2),
$
\end{center}
and their duals. Note that the last two pairs are the symmetric pairs of
the first two. As in previous arguments, it is enough to consider
$\theta_y=t$ or $\pi-t$.

We first handle the triple (in angle notation) 
$(5\pi/12+t/2, t, \pi/12-t/2).$
Note that three of the angles 
in the $\tau_Y$-equation of this point
(all angles rewritten in $[0,\pi/2]$)
are: $5\pi/12+t/2,$ $\pi/12+t/2,$ and $|\pi/12 -3t/2|.$
However,  both $|\pi/12-3t/2|$ and $\pi/12 +t/2$
correspond to angles in $[0,\pi/6).$
If  $|\pi/12-3t/2|$ is non-zero, then the $\tau_Y$-equation of this point
cannot correspond to equations 1-5
since no equation in Theorem \ref{thm:4}
has two angles in $[0,\pi/6).$
If  $|\pi/12-3t/2|=0,$ then
Proposition \ref{prop:5.2} applies. 
For $(5\pi/12+t/2,\pi-t, \pi/12-t/2),$ 
the angles
in the $\tau_Y$-equation of this point
 (all angles rewritten in $[0,\pi/2]$)
are again: $5\pi/12+t/2,$ $|-\pi/12+3t/2|,$ and $\pi/12 +t/2.$

The dual of the pair
$(5\pi/12+t/2, \pi/12-t/2)$ is 
$(7\pi/12-t/2, 11\pi/12+t/2).$
Three of the angles in the $\tau_Y$-equation of
 $(7\pi/12-t/2,t, 11\pi/12+t/2)$ 
are: 
$7\pi/12-t/2,$ $11\pi/12+3t/2,$ and  $11\pi/12-t/2,$
which, when rewritten in $[0,\pi/2],$ become:
$5\pi/12+t/2,$ $|\pi/12-3t/2|,$ and 
 $\pi/12 +t/2.$ These are the same angles as before, thus the same
argument applies.
The same argument also holds for $(7\pi/12-t/2, \pi - t, 11\pi/12+t/2).$

Consider the symmetric pair
$(\pi/12-t/2, 5\pi/12+t/2).$
For the triple $(\pi/12-t/2, t, 5\pi/12+t/2),$
three of the angles in the $\tau_Y$-equation of this point 
are: $\pi/12-t/2,$ $5\pi/12+3t/2,$ and $5\pi/12 -t/2.$
Note that the angle $5\pi/12+3t/2,$ rewritten in $[0,\pi/2],$
is:  
$$
\left\{
\begin{array}{ll}
  5\pi/12+3t/2    & t \in (0,\pi/18] \\
  7\pi/12-3t/2  & t\in [\pi/18,\pi/6).\\
\end{array}
\right.
$$
This angle is in $(\pi/3,\pi/2].$
Further note that the angle  $5\pi/12 -t/2$
is in $(\pi/3,\pi/2).$
A calculation rules out equation 1.
Now of the remaining equations, the angle
$\pi/12-t/2$ can only correspond to the angle $\pi/15$
in equation 4.
However, equation 4 does not have any angles inside
$(\pi/3,\pi/2).$
The same argument holds for the triple
$(\pi/12-t/2, \pi-t, 5\pi/12+t/2).$

The dual of the pair $(\pi/12-t/2, 5\pi/12+t/2)$
is $(11\pi/12+t/2, 7\pi/12 - t/2)$
and the $\tau_Y$-equation of the triples
$(11\pi/12+t/2, t,  7\pi/12 - t/2)$ and 
$(11\pi/12+t/2, \pi-t, 7\pi/12 - t/2)$
have the same three angles (in $[0,\pi/2]$) as above.

A similar argument holds for
for the triples associated with the pairs
$(7\pi/12-t/2, \pi/12-t/2)$ and
$(\pi/12-t/2, 7\pi/12-t/2).$

Case 5:
 For $a=\pi/3-t,$ and $b= \pi/2,$
the possibilities are: 
\begin{center}
$(5\pi/12-t/2, \pi/12+t/2), (7\pi/12+t/2, \pi/12+t/2),$
\end{center}
\begin{center}
$(\pi/12+t/2,5\pi/12-t/2),(\pi/12+t/2,7\pi/12+t/2),$
\end{center}
and their duals. Note that the two last pairs are the symmetric pairs of
the first two. As in previous arguments,
it is enough to consider
$\theta_y = t$ or $\pi -t.$ 

We first handle the triple (in angle notation) 
$(5\pi/12-t/2, t, \pi/12+t/2).$
Note that three of the angles 
in the $\tau_Y$-equation of this point 
are: $5\pi/12-t/2,$ $\pi/12+3t/2,$ and $\pi/12 -t/2.$
A calculation shows that equation 1 is ruled out. 
Thus, $\pi/12-t/2$ can only correspond to the angle $\pi/15$
in equation 4,
while $5\pi/12-t/2$ is in $(\pi/3,\pi/2)$ which rules out 
equation 4. The same holds for
the triple $(5\pi/12-t/2, \pi-t, \pi/12+t/2)$ and those associated
with the dual of $(5\pi/12-t/2, \pi/12+t/2).$

We now consider the triples associated with the
symmetric pair $(\pi/12+t/2, 5\pi/12-t/2).$
Consider the triple
$(\pi/12+t/2, t, 5\pi/12-t/2).$
Three of the angles in the $\tau_Y$-equation of this point 
are: $\pi/12+t/2,$ $5\pi/12-3t/2,$ and $5\pi/12 +t/2.$
One can check that equation 1 is ruled out. 
Now
$5\pi/12+t/2$ can only correspond to the angles $7\pi/15$
and $3\pi/7$
in equations 3 and 5.
That is, $t=\pi/10$ or $\pi/42.$ But in either case, the angle
$\pi/12+t/2$ will not appear in any of equations 2-5.
The same holds for
the triple $(\pi/12+t/2, \pi-t, 5\pi/12-t/2)$ and those associated
with the dual of $(\pi/12+t/2, 5\pi/12-t/2).$

Similar arguments hold for triples associated with the pairs
$(7\pi/12+t/2, \pi/12+t/2)$ and
$(\pi/12+t/2, 7\pi/12+t/2).$

Case 6: For $a=t,$ and $b= \pi/2,$
the possibilities are: 
\begin{center}
$(\pi/4+t/2, \pi/4-t/2), (3\pi/4-t/2, \pi/4-t/2),$
\end{center}
\begin{center}
$(\pi/4-t/2,\pi/4+t/2),(\pi/4-t/2,3\pi/4-t/2),$
\end{center}
and their duals. Note that the two last pairs are the symmetric pairs of
the first two.
As before, it is enough to consider
$\theta_y = \pi/3-t$ or $2\pi/ 3 +t.$ 

We first handle the triple (in angle notation) 
$(\pi/4+t/2,\pi/3-t, \pi/4-t/2).$
Note that three of the angles 
in the $\tau_Y$-equation are: $\pi/4+t/2,$ $7\pi/12-3t/2,$ and $\pi/12 -t/2.$
It is clear that equation 1 is ruled out. 
The angle $\pi/12-t/2$ can only correspond to the angle $\pi/15$
in equation 4,
while $7\pi/12-3t/2,$ when rewritten in $[0,\pi/2],$
corresponds to an angle
in $(\pi/3,\pi/2],$ which rules out 
equation 4. 
The same holds for
the triple $(\pi/4+t/2,2\pi/3+t, \pi/4-t/2)$
and those associated
with the dual of $(\pi/4+t/2, \pi/4-t/2).$

Consider the symmetric pair $(\pi/4-t/2, \pi/4+t/2).$
For the triple $(\pi/4-t/2, \pi/3+t, \pi/4+t/2),$
three of the angles (rewritten in $[0\pi/2]$)
in the $\tau_Y$-equation are:
$\pi/4-t/2$, $\pi/12+t/2,$ and $5\pi/12-3t/2.$ 
Note that both 
$\pi/4-t/2$ and $\pi/12+t/2$ are less than $\pi/4,$
leaving only equations 1 and 4. However, a calculation
rules out equation 4, while equation 1 is clearly ruled out.
The same holds for
the triple $(\pi/4-t/2,2\pi/3+t, \pi/4+t/2)$
and those associated
with the dual of $(\pi/4-t/2, \pi/4+t/2).$

Similar arguments hold
for triples associated with the pairs
$(3\pi/4-t/2, \pi/4-t/2)$ and
$(\pi/4-t/2,3\pi/4-t/2).$

For the other three equations, we simply list the solutions 
$(\theta_x,\theta_y,\theta_z)$  
in their respective category.  For simplicity, we do not list
the solutions correspond to the $\Pin$ representations and
those with one global coordinate equal to zero.
Note that the complete set of solutions 
include the dual, symmetric, and dual-symmetric solutions 
in the first and third coordinates to those listed below.

$$
\cos(\pi/5) - \cos(2\pi/5) - \cos(\pi/3)  =  0
$$
\begin{enumerate}
\item ($k \ge 2$) $(\pi/5, \pi/3, 3\pi/5)$, $(\pi/5, 2\pi/3, 2\pi/5)$.
\item (triples appearing in Table 1) $(\pi/5, 2\pi/3, 3\pi/5)$, $(\pi/5, \pi/3,2\pi/5)$.
\item  (triples with an angle not appearing in Theorem~\ref{thm:4})
$(3\pi/10,\pi/3,\pi/10)$, $(3\pi/10,2\pi/3,\pi/10)$, 
$(7\pi/10,\pi/3,\pi/10)$, $(7\pi/10,2\pi/3,\pi/10)$,
$(7\pi/30,2\pi/5,13\pi/30)$, $(7\pi/30,3\pi/5,13\pi/30)$,
$(7\pi/30,2\pi/5,17\pi/30)$, $(7\pi/30,3\pi/5,17\pi/30)$,
$(7\pi/20,2\pi/5,3\pi/20)$, $(7\pi/20,3\pi/5,3\pi/20)$,
$(7\pi/20,\pi/3,3\pi/20)$, $(7\pi/20,2\pi/3,3\pi/20)$,
$(3\pi/20,2\pi/5,13\pi/20)$, $(3\pi/20,3\pi/5,13\pi/20)$,
$(3\pi/20,\pi/3,13\pi/20)$, $(3\pi/20,2\pi/3,13\pi/20)$,
$(11\pi/30,\pi/5,\pi/30)$, $(11\pi/30,4\pi/5,\pi/30)$,
$(19\pi/30,\pi/5,\pi/30)$, $(19\pi/30,4\pi/5,\pi/30)$,
$(9\pi/20,\pi/5,\pi/20)$, $(9\pi/20,4\pi/5,\pi/20)$,
$(9\pi/20,\pi/3,\pi/20)$, $(9\pi/20,2\pi/3,\pi/20)$,
$(11\pi/20,\pi/5,\pi/20)$, $(11\pi/20,4\pi/5,\pi/20)$,
$(11\pi/20,\pi/3,\pi/20)$, $(11\pi/20,2\pi/3,\pi/20)$,
$(5\pi/12,\pi/5,\pi/12)$, $(5\pi/12,4\pi/5,\pi/12)$, 
$(5\pi/12,2\pi/5,\pi/12)$, $(5\pi/12,3\pi/5,\pi/12)$,
$(7\pi/12,\pi/5,\pi/12)$, $(7\pi/12,4\pi/5,\pi/12)$, 
$(7\pi/12,2\pi/5,\pi/12)$, $(7\pi/12,3\pi/5,\pi/12)$.

(triples $(x,y,z)$ whose $\tau_Y$-equation 
does
not correspond to any equation in Theorem \ref{thm:4})
$(4\pi/15,2\pi/5,\pi/15)$, $(4\pi/15,3\pi/5,\pi/15)$,
$(11\pi/15,2\pi/5,\pi/15)$, $(11\pi/15,3\pi/5,\pi/15)$,
$(2\pi/15,\pi/5,8\pi/15)$, $(2\pi/15,4\pi/5,8\pi/15)$,
$(2\pi/15,\pi/5,7\pi/15)$, $(2\pi/15,4\pi/5,7\pi/15)$.

\end{enumerate}

Note that the categories above are not  
mutually exclusive.

$$
\cos(\pi/5) - \cos(\pi/15) + \cos(4\pi/15) -\cos(\pi/3) = 0
$$
\begin{enumerate}
\item $(k\ge 2)$ $(\pi/3,4\pi/5,2\pi/5)$, $(3\pi/5,\pi/5,\pi/3).$

\item (triples appearing in Table 1) $(\pi/3,\pi/5,2\pi/5)$, $(\pi/3,\pi/3,2\pi/5)$,
$(\pi/3,2\pi/3,2\pi/5)$, $(3\pi/5,4\pi/5,\pi/3)$,
$(3\pi/5,\pi/3,\pi/3).$

\item  (triples with an angle not appearing in Theorem \ref{thm:4})
$(11\pi/30,4\pi/15,17\pi/30),$ $(11\pi/30,11\pi/15,17\pi/30)$,
$(11\pi/30,\pi/3,17\pi/30),$ $(11\pi/30,2\pi/3,17\pi/30)$,
$(11\pi/30,4\pi/15,13\pi/30),$ $(11\pi/30,11\pi/15,13\pi/30)$,
$(11\pi/30,\pi/3,13\pi/30),$ $(11\pi/30,2\pi/3,13\pi/30)$,
$(7\pi/30,\pi/15,\pi/30)$, $(7\pi/30,14\pi/15,\pi/30)$,
$(7\pi/30,\pi/3,\pi/30)$, $(7\pi/30,2\pi/3,\pi/30)$,
$(23\pi/30,\pi/15,\pi/30)$, $(23\pi/30,14\pi/15,\pi/30)$,
$(23\pi/30,\pi/3,\pi/30)$, $(23\pi/30,2\pi/3,\pi/30)$,
$(7\pi/30,\pi/15,13\pi/30)$, $(7\pi/30,14\pi/15,13\pi/30)$,
$(7\pi/30,4\pi/15,13\pi/30)$, $(7\pi/30,11\pi/15,13\pi/30)$,
$(7\pi/30,\pi/15,17\pi/30)$, $(7\pi/30,14\pi/15,17\pi/30)$,
$(7\pi/30,4\pi/15,17\pi/30)$, $(7\pi/30,11\pi/15,17\pi/30)$,
$(\pi/6,\pi/5,\pi/10)$,  $(\pi/6,4\pi/5,\pi/10)$,
$(\pi/6,\pi/3,\pi/10)$,  $(\pi/6,2\pi/3,\pi/10)$,
$(5\pi/6,\pi/5,\pi/10)$,  $(5\pi/6,4\pi/5,\pi/10)$,
$(5\pi/6,\pi/3,\pi/10)$,  $(5\pi/6,2\pi/3,\pi/10)$,
$(3\pi/10,\pi/5,11\pi/30)$, $(3\pi/10,4\pi/5,11\pi/30)$,
$(3\pi/10,4\pi/15,11\pi/30)$, $(3\pi/10,11\pi/15,11\pi/30)$,
$(3\pi/10,\pi/5,19\pi/30)$, $(3\pi/10,4\pi/5,19\pi/30)$,
$(3\pi/10,4\pi/15,19\pi/30)$, $(3\pi/10,11\pi/15,19\pi/30)$,
$(3\pi/10,\pi/5,\pi/30)$, $(3\pi/10,4\pi/5,\pi/30)$,
$(3\pi/10,\pi/15,\pi/30)$, $(3\pi/10,14\pi/15,\pi/30)$,
$(7\pi/10,\pi/5,\pi/30)$, $(7\pi/10,4\pi/5,\pi/30)$,
$(7\pi/10,\pi/15,\pi/30)$, $(7\pi/10,14\pi/15,\pi/30)$.

(triples $(x,y,z)$ whose $\tau_Y$-equation 
does not correspond to any equation in Theorem \ref{thm:4})
$(2\pi/15,4\pi/15,\pi/15),$ $(2\pi/15,11\pi/15,\pi/15)$,
$(2\pi/15,\pi/3,\pi/15),$ $(2\pi/15,2\pi/3,\pi/15)$,
$(13\pi/15,4\pi/15,\pi/15),$ $(13\pi/15,11\pi/15,\pi/15)$,
$(13\pi/15,\pi/3,\pi/15),$ $(13\pi/15,2\pi/3,\pi/15)$,
$(4\pi/15,\pi/15,7\pi/15)$, $(4\pi/15,14\pi/15,7\pi/15)$,
$(4\pi/15,\pi/3,7\pi/15)$, $(4\pi/15,2\pi/3,7\pi/15)$,
$(4\pi/15,\pi/15,8\pi/15)$, $(4\pi/15,14\pi/15,8\pi/15)$,
$(4\pi/15,\pi/3,8\pi/15)$, $(4\pi/15,2\pi/3,8\pi/15)$,
$(\pi/5,\pi/15,7\pi/15)$, $(\pi/5,14\pi/15,7\pi/15)$,
$(\pi/5,4\pi/15,2\pi/15)$, $(\pi/5,11\pi/15,2\pi/15)$,
$(8\pi/15,\pi/15,\pi/5)$,  $(4\pi/5,11\pi/15, 2\pi/15)$,
$(4\pi/5,4\pi/15,2\pi/15)$,
$(8\pi/15,14\pi/15,\pi/5)$.

(triples with infinite
orbits, by Proposition~\ref{prop:5.2}.)
$(\pi/5,\pi/5,2\pi/15)$, $(\pi/5,4\pi/5, 2\pi/15)$,
$(4\pi/5,\pi/5,2\pi/15)$, $(4\pi/5,4\pi/5, 2\pi/15)$,
$(3\pi/5,2\pi/3,\pi/3)$,
$(4\pi/15,\pi/15,\pi/15)$,
$(4\pi/15,14\pi/15,\pi/15)$,
$(4\pi/15,4\pi/15,\pi/15)$,
 $(4\pi/15,11\pi/15,\pi/15)$,
$(11\pi/15,\pi/15,\pi/15)$, $(11\pi/15,14\pi/15,\pi/15)$,
$(11\pi/15,4\pi/15,\pi/15)$, $(11\pi/15,11\pi/15,\pi/15)$,
$(\pi/5,\pi/5,7\pi/15)$, $(\pi/5,4\pi/5,7\pi/15)$, 
$(8\pi/15,\pi/5,\pi/5)$, $(8\pi/15,4\pi/5,\pi/5)$.
\end{enumerate}

$$
-\cos(2\pi/5)+\cos(2\pi/15)-\cos(7\pi/15)-\cos(\pi/3) = 0
$$
\begin{enumerate}
\item $(k\ge 2)$ $(\pi/5,3\pi/5,\pi/3)$, $(\pi/5,2\pi/3,\pi/3)$,
$(2\pi/3,2\pi/5,\pi/5)$, $(2\pi/3,\pi/3,\pi/5)$.
\item (triples appearing in Table 1) $(\pi/5,2\pi/5,\pi/3)$, $(2\pi/3,2\pi/3,\pi/5)$
$(\pi/5,\pi/3,\pi/3)$, $(2\pi/3,3\pi/5,\pi/5)$.
\item  (triples with an angle not appearing in Theorem \ref{thm:4})
$(7\pi/30,7\pi/15,19\pi/30)$, $(7\pi/30,\pi/3,19\pi/30)$,
$(7\pi/30,8\pi/15,19\pi/30)$, $(7\pi/30,2\pi/3,19\pi/30)$,
$(7\pi/30,7\pi/15,11\pi/30)$, $(7\pi/30,\pi/3,11\pi/30)$,
$(7\pi/30,8\pi/15,11\pi/30)$, $(7\pi/30,2\pi/3,11\pi/30)$,
$(13\pi/30,2\pi/15,\pi/30)$, $(13\pi/30,\pi/3,\pi/30)$,
$(13\pi/30,13\pi/15,\pi/30)$, $(13\pi/30,2\pi/3,\pi/30)$,
$(17\pi/30,2\pi/15,\pi/30)$, $(17\pi/30,\pi/3,\pi/30)$,
$(17\pi/30,13\pi/15,\pi/30)$, $(17\pi/30,2\pi/3,\pi/30)$,
$(11\pi/30,2\pi/15,\pi/30)$, $(11\pi/30,7\pi/15,\pi/30)$, 
$(11\pi/30,13\pi/15,\pi/30)$, $(11\pi/30,8\pi/15,\pi/30)$, 
$(19\pi/30,2\pi/15,\pi/30)$, $(19\pi/30,7\pi/15,\pi/30)$, 
$(19\pi/30,13\pi/15,\pi/30)$, $(19\pi/30,8\pi/15,\pi/30)$, 
$(3\pi/10,2\pi/5,\pi/6)$, $(3\pi/10,\pi/3,\pi/6)$, 
$(3\pi/10,3\pi/5,\pi/6)$, $(3\pi/10,2\pi/3,\pi/6)$, 
$(7\pi/10,2\pi/5,\pi/6)$, $(7\pi/10,\pi/3,\pi/6)$, 
$(7\pi/10,3\pi/5,\pi/6)$, $(7\pi/10,2\pi/3,\pi/6)$, 
$(7\pi/30,2\pi/5,\pi/10)$, $(7\pi/30,7\pi/15,\pi/10)$,
$(7\pi/30,3\pi/5,\pi/10)$, $(7\pi/30,8\pi/15,\pi/10)$,
$(23\pi/30,2\pi/5,\pi/10)$, $(23\pi/30,7\pi/15,\pi/10)$,
$(23\pi/30,3\pi/5,\pi/10)$, $(23\pi/30,8\pi/15,\pi/10)$,
$(\pi/10,2\pi/5,17\pi/30)$, $(\pi/10,2\pi/15,17\pi/30)$,
$(\pi/10,3\pi/5,17\pi/30)$, $(\pi/10,13\pi/15,17\pi/30)$,
$(\pi/10,2\pi/5,13\pi/30)$, $(\pi/10,2\pi/15,13\pi/30)$,
$(\pi/10,3\pi/5,13\pi/30)$, $(\pi/10,13\pi/15,13\pi/30)$.

(triples $(x,y,z)$ whose $\tau_Y$-equation
does not correspond to any equation in Theorem \ref{thm:4})
$(4\pi/15,7\pi/15, 2\pi/15)$, $(4\pi/15,\pi/3,2\pi/15)$,
$(4\pi/15,8\pi/15,2\pi/15)$, $(4\pi/15,2\pi/3,2\pi/15)$,
$(11\pi/15,7\pi/15,2\pi/15)$, $(11\pi/15,\pi/3,2\pi/15)$,
$(11\pi/15,8\pi/15,2\pi/15)$, $(11\pi/15,2\pi/3,2\pi/15)$,
$(\pi/15,2\pi/15,7\pi/15)$, $(\pi/15,\pi/3,7\pi/15)$,
$(\pi/15,13\pi/15,7\pi/15)$, $(\pi/15,2\pi/3,7\pi/15)$,
$(\pi/15,2\pi/15,8\pi/15)$, $(\pi/15,\pi/3,8\pi/15)$,
$(\pi/15,13\pi/15,8\pi/15)$, $(\pi/15,2\pi/3,8\pi/15)$,
$(4\pi/15,7\pi/15,2\pi/5)$, $(4\pi/15,8\pi/15,2\pi/5)$,
$(4\pi/15,7\pi/15,3\pi/5)$, $(4\pi/15,8\pi/15,3\pi/5)$,
$(2\pi/5,2\pi/15,\pi/15)$, $(2\pi/5,13\pi/15,\pi/15)$,
$(3\pi/5,2\pi/15,\pi/15)$, $(3\pi/5,13\pi/15,\pi/15)$,

(triples with infinite
orbits, by Proposition~\ref{prop:5.2}.)
$(4\pi/15,3\pi/5,3\pi/5)$,
$(4\pi/15,2\pi/5,3\pi/5)$,
$(3\pi/5,3\pi/5,\pi/15)$, 
$(3\pi/5,2\pi/5,\pi/15)$, 
$(2\pi/5,3\pi/5,\pi/15)$,
$(2\pi/5,2\pi/5,\pi/15)$,
$(2\pi/15,13\pi/15,8\pi/15)$,
 $(2\pi/15,13\pi/15,7\pi/15)$, 
$(4\pi/15,3\pi/5,2\pi/5)$, 
$(4\pi/15,2\pi/5,2\pi/5)$,
$(2\pi/15,2\pi/15,7\pi/15)$,
$(2\pi/15,2\pi/15,8\pi/15)$,
$(2\pi/15,7\pi/15,7\pi/15)$,
$(2\pi/15,8\pi/15,7\pi/15)$,
 $(2\pi/15,8\pi/15,8\pi/15)$,
 $(2\pi/15,7\pi/15,8\pi/15)$,
\end{enumerate}

\end{document}